\documentclass[12pt]{article}
\usepackage{geometry}                
\usepackage{amsmath, amssymb, amsthm, enumerate}

\newtheorem{theorem}{Theorem}[section]
\newtheorem{lemma}[theorem]{Lemma}
\newtheorem{corollary}[theorem]{Corollary}
\newtheorem{proposition}[theorem]{Proposition}

\newtheorem{conjecture}[theorem]{Conjecture}
\newtheorem{problem}[theorem]{Problem}

\newtheorem{question}[theorem]{Question}

\renewcommand{\epsilon}{\varepsilon}

\newcommand{\defword}{\emph}

\title{Oriented trees and paths in digraphs}

\author{Maya Stein\thanks{The author acknowledges support by FONDECYT Regular Grant 1221905 and by ANID Basal Grant CMM FB210005.}\\ \small University of Chile}

\date{}

\begin{document}

\maketitle


\begin{abstract}
Which conditions ensure that a digraph   contains all oriented paths of some given length, or even a all oriented trees of some given size, as a subgraph? One possible condition could be that the host digraph is a tournament of a certain order. In arbitrary digraphs and  oriented graphs, conditions on the chromatic number, on the edge density, on the minimum outdegree and on the  minimum semidegree  have been proposed. In this survey, we  review the known results, and highlight some open questions in the area.
\end{abstract}

\tableofcontents

\section{Introduction}
A main focus of extremal graph theory is the question how a bound on some invariant of a graph $G$ can ensure that $G$ contains a certain graph $H$ as a subgraph. One of the easiest cases is when $H$ is a tree, or even a path. 
For instance, a theorem by Erd\H os and Gallai~\cite{ErdGall59} states that fro all $n,k\in \mathbb N$, every $n$-vertex graph~$G$ with more than $(k-1)n/2$ edges  contains a path on $k$ edges, and in the  famous Erd\H os-S\'os conjecture~\cite{Erdos64} it is suggested that with this density condition we do not only find the $k$-edge path but in fact every $k$-edge tree in~$G$.  Another example is Dirac's theorem~\cite{dirac}, which  establishes a bound on the minimum degree of a graph $G$ that ensures $G$ contains a spanning cycle (or path), and the Koml\'os--S\'ark\"ozy--Szemer\'edi theorem~\cite{KSS2001} gives a similar bound that ensures $G$ contains any spanning tree of bounded degree.

This survey is centered on the question how results of this flavour translate to digraphs. We will consider sufficient conditions for a digraph or an oriented graph to contain all oriented paths or oriented  trees of some fixed order.

A special type of digraph, which has received much attention, is the tournament. (For all notation see Section~\ref{sec:notation}.) Unlike the complete graph $K_n$, which contains all graphs of the same order $n$, a tournament does not contain all oriented graphs of the same or of smaller  order. For instance, the transitive tournament fails to contain any directed cycle.
 On the other hand, R\'edei's theorem~\cite{redei} establishes  that every tournament contains a spanning directed path. Moreover, every tournament of order greater than $7$ contains every  oriented path of the same order~\cite{HT00b}. 
 
There are many more results and open problems concerning subgraphs of tournaments. 
For containment of oriented trees in tournaments, 
 there is a famous conjecture by  Sumner (see~\cite{RW83}), stating that every $2k$-vertex tournament contains each oriented $k$-edge tree. This conjecture was  solved for large $n$ in~\cite{KMO:ProofSumner}. Newer variants take the number of leaves of the tree or its maximum degree into account.  We will survey these developments in Section~\ref{sec:tourn}.
  
  Another prominent class of digraphs, with respect to containment of oriented trees, are oriented graphs of high chromatic number. 
   The famous Gallai-Hasse-Roy-Vitaver theorem, rediscovered several times, states that there is  a directed path of length $k$ in every orientation of a graph $G$ if and only if the chromatic number of~$G$ is $k+1$. It was conjectured in~\cite{BurrO} that every oriented graph of chromatic number~$k+1$ should in fact  contain {\it any} oriented path of length $k$.
  
  Let us turn from containment of oriented paths to containment of oriented trees. 
    In this context, a well-known conjecture of Burr~\cite{BurrO} states that any oriented graph~$D$ of chromatic number  $2k$ should contain any oriented $k$-edge tree. This would imply Sumner's conjecture.
   The impact  of the chromatic number of an oriented graph on the appearance of  oriented paths and trees in it and similar open questions  are discussed in Section~\ref{sec:chrom}.
    
    Our next focus are conditions on the edge density of the host digraph. That is, we will look for analogues of the Erd\H os-S\'os conjecture or similar statements for digraphs. However, we immediately encounter a huge restriction. Namely, a condition on the edge density alone can only force \defword{antidirected subdigraphs}, that is,~subdigraphs without any vertices of both positive   outdegree and positive indegree.  To see this, consider 
     a complete bipartite graph on two vertex sets $A$ and $B$ and direct all edges from $A$ to $B$. The obtained digraph has  $|A||B|$ edges, and   all its subdigraphs are antidirected.
     
      So, unless we wish to impose extremely high bounds on the number of edges of the host digraph, the only oriented trees that can be forced by a density condition are the  antidirected trees. A beautiful conjecture in this direction appeared in~\cite{ABHLSTR}. It states that every digraph with more than $(k-1)n$ edges contains each antidirected $k$-edge tree. In a certain sense, this would unify the Erd\H os--S\'os conjecture and Burr's conjecture for antidirected trees.  See Section~\ref{sec:dens} for more details on this conjecture and other results.
 
We now turn to conditions on analogues of the minimum degree in digraphs. 
 Conditions on the minimum outdegree that imply the existence of oriented paths have been considered in the literature, but these problems appear to be very difficult. For instance, Thomass\'e conjectured~\cite{havetGarden2, summaryCH} that every oriented graph of minimum outdegree at least $k$ contains a directed path of length $2k$, and this seems wide open (see e.g.~\cite{Subdivi} for related results concerning subdivisions). 

This changes if we add a condition on the minimum indegree.
In other words, 
we will  consider conditions on the \defword{mini\-mum semidegree} of a digraph that can ensure it contains oriented paths or trees of a given order. This will be the topic of Section~\ref{sec:semi}. 

In contrast to Thomass\'e's conjecture~\cite{havetGarden2, summaryCH} mentioned above, it is known that for any even $k\in\mathbb N$, every oriented graph of minimum semidegree at least $k/2$ contains a directed path of length $k$, this is a theorem of Jackson~\cite{jackson}. 
A generalisation of this statement to all oriented paths was recently conjectured in~\cite{tree-survey}: for every $k\in\mathbb N$, every oriented graph of minimum semidegree exceeding $k/2$ should contain any oriented path of length $k$.  The analogous question for hosts that are digraphs instead of oriented graphs is more complicated, because a disjoint union of complete digraphs or order $k$ each has minimum semidegree  $k-1$ but contains no oriented subgraph with $k$ edges. 

We  move on to results and open problems concerning semidegree conditions for finding oriented trees.  The central result here is the recent generalisation~\cite{km} of the 
 Koml\'os--S\'ark\"ozy--Szemer\'edi theorem~\cite{KSS2001} to digraphs. This result assures that every large enough digraph obeying a certain minimum degree condition contains each oriented tree of the same order, if the maximum total degree of the underlying tree is bounded. Another recent result concerns balanced antitrees of bounded maximum degree in oriented graphs~\cite{camila}, where the antitree may be of lower order than the host.

Finally, in Section~\ref{sec:pseu}, we investigate another degree condition, the \defword{minimum pseudo-semidegree}.
For a non-empty digraph $D$, this is the maximum $k$ such that each vertex of $D$ has either outdegree $0$ or outdegree at least $k$, and has either indegree $0$ or indegree at least $k$. 
The motivation for defining the minimum pseudo-semidegree is that, as we will see in 
Section~\ref{sec:pseu}, any digraph of large edge-density has a subgraph of large minimum pseudo-semidegree, and this is not true if we consider the minimum semidegree instead.
 Moreover,  if we wish to find antidirected subgraphs in a digraph, bounds on the minimum semidegree and bounds on the minimum pseudo-semidegree have almost the same effect (see Lemma~\ref{antiiiilemma}). Consequently, some of the results from Section~\ref{sec:semi} already hold if we consider the minimum pseudo-semidegree instead of the minimum semidegree.

 \section{Notation}\label{sec:notation}
 For easy reference, we gather all definitions here, except for a few that are only used once in specific settings. All of the definitions are standard.
 
 \paragraph{Digraphs and oriented graphs.}
 A \defword{digraph} $D$ consists of a set $V(D)$ of \defword{vertices}, and a set  $E(D)$ of \defword{directed edges} (or \defword{edges} for short) which correspond to pairs of distinct vertices. We write $ab$ for the edge $(a,b)$, that is, the edge from $a$ to $b$, and note that a digraph can have both edges $ab$ and $ba$. 
 
  A digraph is called \defword{symmetric} if $ab\in E(D)$ implies that $ba\in E(D)$. So, there is a one-to-one correspondence between symmetric digraphs and graphs.
  
 An \defword{oriented graph} is a digraph that has at most one edge between any pair of vertices. Thus, an oriented graph $D$ corresponds to an \defword{orientation} of an undirected graph $G$. We also call $G$ the \defword{underlying graph} of $D$.

\paragraph{Degrees.}
In a digraph $D$, the \defword{outdegree} $d^+(v)$ of a vertex $v$ is the number of edges of~$D$ that are directed from $v$ to another vertex. The \defword{indegree} $d^-(v)$ of a vertex is defined analogously.  

 An oriented graph is called \defword{$\ell$-regular} if each of its vertices has outdegree $\ell$ and indegree $\ell$.
 A {\it regular} tournament is an $n$-vertex tournament which is $(n-1)/2$-regular, for any odd $n\in\mathbb N$. Using induction, it is easy to see that regular tournaments on $n$ vertices exist for each odd $n\in\mathbb N$.

The
 \defword{minimum outdegree} $\delta^+(D)$ of a digraph $D$ is defined as the minimum over all the outdegrees of the vertices, and the \defword{minimum indegree} $d^-(D)$  is defined ana\-logously.
The \defword{total degree} of a vertex is the sum of its in- and its outdegree. 
The {\it maximum total degree} $\Delta(D)$ of an oriented graph $D$ is defined as the maximum over the total degrees of the vertices, that is, the maximum degree of the underlying undirected graph.
Another   widely used degree notion is the
 \defword{minimum semidegree} $\delta^0(D)$ of a digraph $D$, which is defined as the minimum over all the in- and all the out-degrees of the vertices. 

%
%

\paragraph{Orientations.}
In graphs, if a path has $k$ edges, we say it is a \defword{$k$-edge path} or that it has \defword{length} $k$, and use the same terminology for oriented paths. For $k\in \mathbb N $, the  \defword{directed $k$-edge path} (or \defword{directed path on $k$ edges})  is the  $k$-edge path having all its edges oriented in the same direction. A \defword{directed cycle} is a cycle with all edges oriented in the same direction. 

 An {\it out-arborescence} is an oriented tree having all its edges directed away from a specific vertex. More precisely, an arborescence is an orientation of a tree $T$ with a root $r$ such that an edge $xy$ is directed from $x$ to $y$ if and only if $x$ lies on the unique path from $y$ to $r$. An {\it in-arborescence} is defined analogously. Note that a directed path qualifies as both an out- and an in-arborescence
 
 Another example of an out-arborescence is  the $k$-edge star with all its edges directed outwards, that is, away from its centre. This oriented star is denoted by  $K_{1\rightarrow k}$. Similarly,  $K_{1\leftarrow k}$ is the  $k$-edge star with all its edges directed  towards its centre.
 
 An oriented graph $D$ is called \defword{antidirected} if no vertex has both positive outdegree and positive indegree. 
 So, any antidirected graph
  has a bipartition $(A,B)$ which induces the directions on the edges:  all edges are directed from $A$ to~$B$.  
 Note that if $T$ is an oriented tree, this bipartition coincides with the natural bipartition of the underlying tree.
   Finally, a tree, or an oriented tree, is called {\it balanced}  if its bipartiton classes have the same size.

  \paragraph{Complete digraphs and tournaments.}
  A digraph is called \defword{complete} if it contains all possible edges, that is, between each pair of vertices, both edges are present. 
  
  An oriented graph is called a \defword{tournament} if it contains exactly one edge between each pair of vertices. 
   The \defword{transitive tournament} is the tournament on vertices $v_1, \ldots , v_n$, with all edges directed towards the endpoint with the larger index. Alternatively, one can define a transitive tournament as a tournament that contains no directed cycles.
   
      \paragraph{Other notation.}
   
   The \defword{chromatic number} of an oriented graph is the chromatic number of the underlying graph.
 The chromatic number of an arbitrary digraph $D$ is defined analogously, considering the graph obtained by identifying multiple edges of the underlying multigraph of $D$.
 
 A vertex or a set of a digraph $D$ is said to \defword{dominate} another set or vertex if all edges from the former to the latter are present.  A digraph is \defword{strong} if for each pair $x,y$ of distinct vertices there is a directed path starting in $x$ and ending in $y$.
 
  A {\it blow-up} of a digraph $D$ is obtained by replacing each vertex with an independent set of vertices, and adding all edges from such a set $X$ to a set $Y$, if  $X$ and $Y$ originated from vertices $x$ and $y$ belonging to an edge ${xy}$ of $D$. We call a blow-up $D'$ of $D$ an \defword{$\ell$-blow-up of $D$} if each original vertex from $D$ was replaced with exactly $\ell$ vertices in $D'$.

 \section{Tournaments}\label{sec:tourn}
 
\subsection{Paths in tournaments}
The first results on oriented trees in tournaments concern oriented paths, and more specifically, directed paths. %
 One of the earliest results in this respect is R\'edei's theorem~\cite{redei} from~1934, on spanning directed paths in tournaments. 
 \begin{theorem}[R\'edei 1934~\cite{redei}]\label{redei}
 For each $k\in\mathbb N$, every tournament on $k+1$ vertices contains the directed $k$-edge path.
 \end{theorem}
 The proof of this theorem is simple: Take a longest directed path $P=v_1v_2\ldots v_\ell$, and assume there is a vertex $x\notin V(P)$. Let $i$ be the first index such that $xv_i$ is an edge, and if no such index exists, set $i:=\ell +1$. Then we can insert $x$ between $v_{i-1}$ and $v_i$ (or before $v_i$, or after $v_{i-1}$ if $i=1$ or $i=\ell +1$) to obtain a longer directed path than $P$, a contradiction.
 
 Subsequently, containment of other spanning oriented paths was studied. Gr\"un\-baum~\cite{gruenbaum} showed in 1971 that   each  antidirected $k$-edge path is contained in each tournament on $k+1$ vertices, with three exceptions: 
 the regular tournaments on $3$ and $5$ vertices, and the Paley tournament on $7$ vertices. Let us call the set of these three tournaments $\mathcal T^*$.
In 1972,  Rosenfeld~\cite{rosenfeld} gave an easier proof of Gr\"unbaum's result and showed that moreover, any vertex of the tournament can be chosen as the starting point of the path. 

Further, Rosenfeld~\cite{rosenfeld} conjectured that there is a $n>7$ such that every tournament on at least $k+1\ge n$ vertices contains  all oriented  paths on $k$ edges. Progress towards this conjecture was achieved in~\cite{alspRos, forc, RW83, straight, zhang}, and finally, in 1986, the conjecture was settled by Thomason~\cite{tho86}. He showed that $n=2^{128}$ suffices, while expressing the belief that $n=8$ should be the correct number.    In 2000, Havet and Thomass\'e~\cite{HT00b} confirmed this, and showed that the only exceptional tournaments not containing all oriented paths of the corresponding length are the ones found by Gr\"unbaum in 1971:

 \begin{theorem}[Havet and Thomass\'e 2000~\cite{HT00b}]\label{thmHT}
  For each $k\in\mathbb N$, every tournament $T$ on $k+1$ vertices contains each oriented $k$-edge path $P$, unless $T\in\mathcal T^*$ and $P$ is antidirected.
 \end{theorem}
 
 A shorter proof of Theorem~\ref{thmHT} was recently given by Hanna~\cite{hanna}.
 In~\cite{HT00b}, Havet and Thomass\'e also obtain the additional result that for each vertex $v$ of the tournament $T$ and each oriented $k$-edge path $P$, there is a copy of $P$ in $T$ either starting or ending at~$v$, as long as $|T|\ge 9$.

There are some related questions which have been studied in the literature. To mention a few, the first result  on the number of directed Hamilton paths in a tournament is from 1943 and due to Szele~\cite{szele43}. Szele's lower bound is completemented by an upper bound found by Alon~\cite{alon-number} in 1990. See e.g.~\cite{ElSahGhaHan}
for some more recent results.
In a slightly different direction, Thomassen~\cite{tho82}  characterised in 1982 the tournaments that do not have two directed edge-disjoint spanning paths. Thomassen~\cite{tho80,tho82, tho84} also showed several results on spanning directed paths starting or ending at specified vertices, or through specified edges.
See the survey~\cite{JBJ-GG} for more directions.

\subsection{Trees in tournaments}\label{sec:treetour}
 
 If instead of any oriented path, we wish to guarantee any oriented tree of a certain size in a tournament, it turns out that we need the tournament  to be substantially larger than the tree. 
 In order to see this, let us  consider $K_{1\rightarrow k}$, the $k$-edge star with all its edges directed outwards. For $n<2k$, there are $n$-vertex tournaments of maximum outdegree at most $k-1$. Clearly, these do not contain $K_{1\rightarrow k}$ as a subgraph. So, if we aim for the statement  that all tournaments on $f(k)$ vertices contain all oriented $k$-edge trees, we need to choose $f(k)\ge 2k$.
 
 It is easy to see that $f(k)=2k$ is sufficient for finding $K_{1\rightarrow k}$, for any $k\in\mathbb N^+$, as every tournament on at least $2k$ vertices necessarily has a vertex of outdegree at least~$k$. 
 Clearly, any tournament $T$ on $2k$ vertices also contains the other antidirected star, $K_{1\leftarrow k}$. In fact, it contains any oriented $k$-edge star. Indeed, say 
$S$ is a $k$-edge oriented star with  $\ell$  `out-leaves' (vertices of in-degree $1$).
If $T$ fails to contain $S$, then each vertex of $T$ either has outdegree less than $\ell$ or in-degree less than $k-\ell$ (in cannot fulfill both these conditions, as its total degree is $2k-1$). So $V(T)$ splits into the two sets $V_1=\{v\in V(T):d^+(v)<\ell\}$ and $V_2=\{v\in V(T):d^-(v)<k-\ell\}$. As $V_1$ has $|V_1|(|V_1|-1)/2$ edges, and therefore some $v\in V_1$ sends edges to at least $\lceil(|V_1|-1)/2\rceil$ vertices of~$V_1$, we know that $\lceil(|V_1|-1)/2\rceil<\ell$ or $|V_1|=0=\ell$. It follows that $|V_1|\le 2\ell -1$ or $V_1=\emptyset$. Similarly, we see that $|V_2|\le 2k-2\ell-1$ or $V_2=\emptyset$. So $2k=|V_1|+|V_2|< 2k$, a contradiction. Note that if $\ell\neq0\neq k-\ell$, then this argument actually shows that any tournament on $2k-1$ vertices contains the $k$-edge oriented star with  $\ell$  out-leaves.
 
Sumner conjectured in 1971  that the same bound on the order of the host tournament would be sufficient to find  any oriented tree having $k$ edges:
\begin{conjecture}[Sumner 1971, see~\cite{RW83}]\label{sumner}
Every tournament on $2k$ vertices contains every oriented $k$-edge tree.
\end{conjecture}

Conjecture~\ref{sumner} has also been called `Sumner's universal tournament conjecture' in the literature. 
Much of the early progress towards Conjecture~\ref{sumner} consists of proofs for specific classes of oriented trees, or by replacing the bound $2k$ from  Conjecture~\ref{sumner} with a larger number (see for instance~\cite{chung, DrossHavet, elsahili, HT91,  Havet02, HT00a, worm}). 

The first result showing that any tournament on $O(k)$ vertices contains every oriented $k$-edge tree is due to H\"aggkvist and Thomason~\cite{HT91}. The
 currently best  bound on the order of the host tournament was established in~2021  by Dross and Havet~\cite{DrossHavet}: They proved that for each $k\in\mathbb N$, every tournament on  $\lceil 2.625 k- 2.9375\rceil$ vertices contains each oriented $k$-edge tree. Furthermore, in 2000, Havet and Thomass\'e~\cite{HT00a} showed that Conjecture~\ref{sumner} holds for all trees that are arborescences, which could be viewed as an analogue of R\'edei's theorem for trees.

In 2011, K\"uhn, Mycroft and Osthus proved an approximate version~\cite{KMO:ApproximateSumner} of Sumner's conjecture (Conjecture~\ref{sumner}), 
and shortly afterwards, the same authors confirmed the conjecture  for large $k$:
\begin{theorem}[K\"uhn, Mycroft and Osthus 2011~\cite{KMO:ProofSumner}]\label{KMOsumner}
There is a $k_0\in\mathbb N$ such that for every $k\ge k_0$, every tournament on $2k$ vertices contains every oriented $k$-edge tree.
\end{theorem}

 Interestingly, for most trees, the host tournament can be much smaller. 
If we wish to find a random oriented tree in a tournament, we have the following result, which was conjectured in 1988 by Bender and Wormald~\cite{BenWor} and was proved thirty years later by Mycroft and Naia~\cite{MyNa18}.
 \begin{theorem}[Mycroft and Naia 2018~\cite{MyNa18}]
 Let $T$ be chosen uniformly at random from the set of all labelled oriented trees with $k$ edges. Then asymptotically almost surely every tournament on $k+1$ vertices contains $T$.
 \end{theorem}
 
 So there is a large class of trees which allow for a better bound on the size of the host tournament in Conjecture~\ref{sumner}. 
Since the known extremal examples for Conjecture~\ref{sumner} are the antidirected stars, and in view of Theorems~\ref{redei} and~\ref{thmHT}, which show that Conjecture~\ref{sumner} is far from tight for oriented paths, it seems natural to suspect that  trees which are very different from stars may be good candidates for such  trees. One natural possibility would be to consider a restriction  on the maximum degree of the underlying tree $T$.
 And indeed, it turns out that for oriented trees of bounded maximum degree, the size of the host tournament can be lowered almost to $k$. 
 
 The first results in this direction are for a type of trees that has been called `claws' in the literature. A \defword{claw} is a collection of directed paths whose first vertices have been identified. Improving previous results by Saks and S\'os~\cite{saksSos} and by Lu~\cite{lu, lu2}, it was shown by Lu, Wang and Wong~\cite{lww} in 1998 that every tournament on $k+1$ vertices contains each $k$-edge claw of maximum degree at most $19(k+1)/50$, and that, on the other hand, there is a family of $k$-edge claws with maximum degree converging to $11(k+1)/23$ which are not contained in all tournaments on $k+1$ vertices. 
 
Arbitrary trees of bounded maximum degree were first studied in~\cite{KMO:ApproximateSumner}. 
 K\"uhn, Mycroft and Osthus~\cite{KMO:ApproximateSumner} proved in  2011  that if $k$ is sufficiently large, then every tournament on~$(1+o(1))k$ vertices contains every orientation of every $k$-edge tree~$T$ with  $\Delta(T)\le\Delta$, where $\Delta$ is some constant. 
Mycroft and Naia~\cite{MyNa18} improved this bound to $\Delta(T)\le (\log k)^{\Delta'}$, where $\Delta'$ is a  constant.
 Very recently, Benford and Montgomery gave a linear bound on the maximum degree of the tree.
 
  \begin{theorem}[Benford and Montgomery~\cite{BenMon2}]\label{thm:BenMon2}
For each $\varepsilon$, there is a $k_0\in\mathbb N$ and a constant $c$ such that for every $k\ge k_0$, every tournament on $(1+\varepsilon)k$ vertices contains every orientation of every $k$-edge tree $T$ with $\Delta(T)\le ck$.
\end{theorem}
 
 The following related question first appears explicitly in~\cite{MyNa18}:
\begin{question}[Mycroft and Naia~\cite{MyNa18}]\label{q:myna}
Is  there  a function  $g$ such that every tournament on $k+g(\Delta)$ vertices contains every orientation of every $k$-edge tree $T$ with $\Delta(T)\le \Delta$?
\end{question}

Mycroft and Naia~\cite{MyNa18} also ask whether one can choose $g(\Delta)=2\Delta -4$ if~$k$ is much larger than $\Delta$. 
Earlier examples by Allen and Cooley (see~\cite{KMO:ApproximateSumner})  show  that  this would be best possible: Let $h_1, h_2\in\mathbb N$. For this, consider the oriented tree $T$ obtained 
by taking a directed path $P=v_0v_1\ldots v_{h_2}$, and adding $2h_1$ new vertices, of which $h_1$ vertices send an edge to $v_0$, and $h_1$ vertices receive an edge from $v_{h_1}$. Note that~$T$ has $k:=2h_1+h_2$ edges and $\Delta(T)=h_1+1$. 
Now, take the disjoint union of three  tournaments $H_1$, $H_2$ and $H_3$, where $H_1$ and $H_3$ are $(h_1-1)$-regular tournaments  on $2h_1-1$ vertices, and $H_2$ is any tournament on $h_2$ vertices. For $1\le i<j\le 3$, add all edges between $H_i$ and $H_j$, orienting them from $H_i$ to $H_j$. It is not hard to see that the resulting tournament $H$ does not contain $T$, and $H$ has $4h_1+h_2-3=k+2\Delta(T)-5$ vertices. 

 Looking for conditions that make a tree resemble more a path than a star, one could, 
instead of restricting the maximum degree of $T$,  impose a condition on the number of leaves of $T$. In this very natural direction, 
H\"aggkvist and Thomason~\cite{HT91} showed in 1991 that there is a function $f$ (which is exponential  in $\ell^3$) such that  every tournament on $k+f(\ell)$ contains each oriented tree with $k$ edges and at most~$\ell$ leaves. 
Havet and Thomass\'e proposed in  1996 that it should be possible to drop the size of the tournament from Sumner's conjecture to $k+\ell$ if we are looking for an oriented  tree  having at most $\ell$ leaves. In other words, they proposed the following generalisation of Conjecture~\ref{sumner}.

\begin{conjecture}[Havet and Thomass\'e 1996, see~\cite{Havet03}]\label{conj:HT}
 Let $T$ be an oriented $k$-edge tree
 with $\ell$ leaves. Then every tournament on $k +\ell$ vertices
contains a copy of $T$.
\end{conjecture}

In 2021, Dross and Havet~\cite{DrossHavet} showed that there is a function $f$ which is quadratic in $\ell$ such that every tournament on $k+f(\ell)$ vertices contains each oriented tree with $k$ edges and at most $\ell$ leaves. 
Recently, Benford and Montgomery~\cite{BenMon} proved a linear bound on the number of `extra' vertices in the tournament.

 \begin{theorem}[Benford and Montgomery 2022~\cite{BenMon}]
There is some $C$
 such that every 
$(k+C\ell)$-vertex tournament contains a copy of any $k$-edge oriented tree with~$\ell$ leaves.
 \end{theorem}
 
 Let us turn back to Conjecture~\ref{conj:HT}.
Clearly, for $K_{1\rightarrow k}$ and $K_{1\leftarrow k}$, the bound from Conjecture~\ref{conj:HT} coincides with the bound from Conjecture~\ref{sumner}, and for all other trees, it is lower. 
Moreover, the bound from  Conjecture~\ref{conj:HT} is not tight for small $\ell$: As we saw above in Theorem~\ref{thmHT}, if $T$ is a path, i.e.~a tree with $\ell=2$ leaves, and if $k\ge 7$, then already any tournament on at least $k+1=k+\ell-1$ vertices contains $T$, which is one smaller than required by Conjecture~\ref{conj:HT}. Ceroi and Havet~\cite{CeroiHavet} showed that  for any tree $T$ with $\ell=3$ leaves, if $k\ge 4$, then every tournament on $k+2=k+\ell-1$ vertices contains~$T$. 

Dross and Havet~\cite{DrossHavet} conjectured that the same holds for all trees with few leaves, i.e.~they conjectured that the bound of $k+\ell-1$ is sufficient whenever $k$, the number of edges of  the tree is of sufficiently larger order than $\ell$, its number of leaves. 
This conjecture was confirmed recently by Benford and Montgomery~\cite{BenMon}:
 \begin{theorem}[Benford and Montgomery 2022~\cite{BenMon}]\label{BenMon2}
For every $\ell$ there is a~$k_0$ such that for every $k\ge k_0$, every tournament on $k+\ell-1$ vertices contains each oriented $k$-edge tree with~$\ell$ leaves.
 \end{theorem}
 
 One can understand the number  $k_0$ in Theorem~\ref{BenMon2} being a function of $\ell$. The authors of~\cite{BenMon} state that their $k_0$ can be chosen as $\ell^{O(\ell)}$, but also express their belief that the theorem may stay correct with $k_0$ being of  order $O(\ell)$.


\section{Chromatic number}\label{sec:chrom}
\subsection{Chromatic number and oriented paths}
Recall that the {chromatic number} of an oriented graph is the chromatic number of the underlying undirected  graph. So the chromatic number of a tournament equals its number of vertices. This motivates the question whether the results from Section~\ref{sec:tourn} do not only apply to tournaments, but more generally  to digraphs or oriented graphs of large chromatic number.

An early and quite famous result in this direction is the  Gallai-Hasse-Roy-Vitaver  theorem  (see, e.g.~\cite{Bang-Jensen2018}), discovered several times during the 1960's, which we state next. The theorem actually is a duality result.
 \begin{theorem}[Gallai-Hasse-Roy-Vitaver  theorem, 1960's]\label{thmGHRV}
 Let $G$ be a  graph. Then $G$ has chromatic number $k+1$ if and only if in 
every orientation of  $G$ there is a directed $k$-edge path. 
 \end{theorem}
 We are interested in the following 
 immediate corollary of Theorem~\ref{thmGHRV}, which can be seen as a generalisation of R\'edei's theorem.
\begin{corollary}~\label{coroGHRV}
Every $(k+1)$-chromatic oriented graph contains the  directed path with $k$ edges. 
\end{corollary}

For large $k$, a version of Corollary~\ref{coroGHRV}, replacing directed paths with oriented paths hase been proposed by Burr in~\cite{BurrO} (see also~\cite{chung}):
\begin{conjecture}[Burr 1980]\label{conj:chropath}
 Every $(k+1)$-chromatic oriented graph contains each  oriented  path  with $k$ edges.  
\end{conjecture}

There is some support for this conjecture, apart from it being true for tournaments. El Sahili~\cite{elsah} proved in 2004 that Conjecture~\ref{conj:chropath} holds for $k=3$ and, if the path is antidirected, for $k=4$. Addario-Berry, Havet, and Thomass\'e~\cite{ABHT} proved in 2007 that the conjecture holds for every oriented path with $k \ge 3$ edges that changes direction at most once.

\subsection{Chromatic number and oriented trees}

Let us now turn to oriented trees in digraphs of large chromatic number. 
A variant  of Conjecture~\ref{conj:chropath} replacing   oriented paths with oriented stars is not true (we already saw this for tournaments in Section~\ref{sec:treetour}). However,  as in Sumner's conjecture, it seems natural to consider larger host tournaments.
In this spirit, Burr~\cite{BurrO} suggested the following conjecture, which, if true, would imply Sumner's conjecture.

\begin{conjecture}[Burr 1980~\cite{BurrO}]\label{conj:BurrOrient}
 Every $2k$-chromatic digraph contains  each  oriented $k$-edge tree.
\end{conjecture}

Conjecture~\ref{conj:BurrOrient} is only known for some specific classes of oriented paths (see above) and  for all oriented stars~\cite{naiaThesis}. 
Burr~\cite{BurrO} showed that Conjecture~\ref{conj:BurrOrient} is true if we replace $2k$ with $k^2$. Addario-Berry, Havet, Linhares Sales, Thomass\'e and Reed~\cite{ABHLSTR} improved this bound, roughly by a factor of $1/2$. Naia~\cite{Naia} gives a  bound  which is better for all oriented graphs with large chromatic number, by showing that every $((k+2)\log_2 n)$-chromatic oriented graph contains  each  oriented $k$-edge tree.

A generalisation of Conjecture~\ref{conj:HT} in the spirit of Conjecture~\ref{conj:BurrOrient} has been proposed by Havet~\cite{havetGarden} in 2013, and appears in~\cite{DrossHavet}: 
\begin{conjecture}[Havet 2013~\cite{DrossHavet, havetGarden}]\label{conjHavet}
Every $(k+\ell)$-chromatic digraph contains each  oriented $k$-edge tree having $\ell$ leaves. 
\end{conjecture}

In view of the results for tournaments exhibited in the previous section, it seems natural to ask whether the maximum degree of the oriented tree could have a similar influence on the bound on the chromatic number.
In analogy to Theorem~\ref{thm:BenMon2} and Question~\ref{q:myna}, it seems natural to ask the following.

\begin{question}
For each $\varepsilon$, are there is a $k_0 $ and   $c$ such that for every $k\ge k_0$, every $(1+\varepsilon)k$-chromatic digraph   contains every orientation of every $k$-edge tree $T$ with $\Delta(T)\le ck$?
\end{question}

\begin{question}
Is  there  a function  $g$ such that every $(k+g(\Delta))$-chromatic digraph contains every orientation of every $k$-edge tree $T$ with $\Delta(T)\le \Delta$?
\end{question}

We add that Naia~\cite{Naia} explicitly conjectures that 
for every oriented $k$-edge tree~$T$, the minimum~$n$ such that every tournament of order~$n$ contains $T$ coincides with the minimum $n$ such that every $n$-chromatic oriented graph contains~$T$. 
We close the section by referring the reader to \cite{digraph-like} for a survey of more open problems related to tournaments and containment of trees and other subdigraphs.

 \section{Edge density}\label{sec:dens}
 
For graphs, the influence of the edge density of a graph $G$ on the appearance of paths and other trees in $G$  has been much studied. The first result in this direction is a theorem by Erd\H os and Gallai from 1959~\cite{ErdGall59}, which states that for each $k\in \mathbb N$, each graph of average degree greater than $k-1$ contains the path on $k$ edges. Of course, such a graph also contains the $k$-edge star.

In 1963, Erd\H os and S\'os (see~\cite{Erdos64}) suggested that the same should hold for each tree: They conjectured that for each $k\in \mathbb N$, each graph of average degree greater than $k-1$ contains each tree with $k$ edges. The Erd\H os--S\'os conjecture has been shown to be true for large dense graphs: in an sharp version for large trees of bounded maximum degree~\cite{BPS3}, and very recently, in  an approximate version for  all large trees~\cite{newDiana}. See~\cite{tree-survey} for an overview of this conjecture.

Note that the condition on the average degree in the Erd\H os--Gallai theorem and in the Erd\H os--S\'os conjecture is best possible as a complete graph on $k$ vertices fails to contain a $k$-edge path or tree. Also, it is clear that the condition on the average degree can be reformulated as a condition on the edge density of the host graph $G$, namely the condition that $|E(G)|/|V(G)|>(k-1)/2$.


If one wishes to extend these statments  to digraphs, the first obstacle one encounters is the fact that for all $n\in \mathbb N$ there are oriented graphs $D$ on $2n$ vertices and of  edge density $|E(D)|/|V(D)|=n/2$ that are antidirected, and thus only have antidirected subgraphs. These are the oriented graphs obtained from the complete bipartite graph $K_{n,n}$, for any $n\in \mathbb N$, with all edges directed from the first partition class to the second. So the only oriented $k$-edge  trees that can be guaranteed by imposing a lower bound $f(k)$ on the edge density of the host digraph  are the antidirected trees. (This observation is attributed in~\cite{BurrO} to de Bruijn, and the  explicit example above appears in \cite{ABHLSTR}.)

Thus motivated, 
it becomes natural to ask which edge density is sufficient to guarantee that a digraph contains all antidirected trees of a certain size.
In this direction,
Graham~\cite{graham} confirmed in 1970 a conjecture he attributes to Erd\H os: for every antidirected tree $T$ there is a constant $c_T$ such that every sufficiently large directed graph $D$ on $n$ vertices and with at least $c_Tn$ edges contains~$T$. 
In 1982, Burr~\cite{BurrO} gave an improvement of Graham's result:

 \begin{theorem}[Burr~\cite{BurrO}]\label{thm:4kn}
 Every $n$-vertex digraph~$D$ with more than $4kn$ edges contains each antidirected tree~$T$ on $k$ edges. 
 \end{theorem}
 
Burr obtains Theorem~\ref{thm:4kn} by taking a maximal cut, which contains at least half of the edges of $D$, and then omitting all edges from one of the possible directions in this cut. This gives an oriented bipartite subgraph $D'=(A,B)$  of~$D$ which has more than $kn$ edges, which are all directed in the same way, say from $A$ to $B$. Then the average total degree of the vertices of $D'$ exceeds $2k$. A standard argument yields a subgraph $D''$ of $D'$ such that $D''$ has minimum total degree at least $k$. So in $D''$, each vertex from $A$ has outdegree at least $k$, and each vertex from $B$ has indegree at least $k$. Now, one can embed the $k$-edge tree $T$ greedily into $D''$, thus completing the proof of Theorem~\ref{thm:4kn}. 

Burr writes in~\cite{BurrO} that the bound $4kn$ on the number of edges can `almost certainly be made rather smaller'. He also provides an example of an oriented graph on $n$ vertices wit $(k-1)n$ edges that fails to contain $K_{1\rightarrow k}$, the $k$-edge star with all edges directed outwards. Namely, $K_{1\rightarrow k}$ is not contained in  the complete bipartite graph $K_{2k-2, 2k-2}$ with half of the edges  oriented in either direction in an appropriate way. Note that the graph $K_{2k-2, 2k-2}$ has $n=4k-4$ vertices and  $(2k-2)^2=(k-1)n$ edges. 

In the previous example, instead of $K_{2k-2, 2k-2}$, one could also consider any $(k-1)$-regular tournament on $n$ vertices or the $(k-1)$-blow-up of a directed cycle on $n$ vertices  which each have $(k-1)n$ edges but do not contain $K_{1\rightarrow k}$. Moreover, the complete di\-graph on~$n=k$ vertices has $(k-1)n$ edges and does not contain {\it any} oriented $k$-edge tree. 

In 2013, 
Addario-Berry, Havet, Linhares Sales, Reed and Thomassé~\cite{ABHLSTR}   conjectured   that the bound $(k-1)n$, which was shown to be necessary by Burr, is indeed the correct bound.

\begin{conjecture}[Addario-Berry et al. 2013~\cite{ABHLSTR}]\label{conj:antidir}
 Every  digraph $D$ with more than $(k-1)|V(D)|$ edges contains each  antidirected $k$-edge tree.
\end{conjecture}

It is observed in~\cite{ABHLSTR} is that Conjecture~\ref{conj:antidir} implies Conjecture~\ref{conj:BurrOrient}  for antidirected trees, because  every $2k$-chromatic graph has a subgraph $H$ of minimum degree at least $2k-1$, and thus  $H$ has more than $(k-1)|V(H)|$ edges.
Moreover,  
if
we restrict  Conjecture~\ref{conj:antidir} to symmetric digraphs, then it becomes equivalent to the Erd\H os--S\'os conjecture mentioned at the beginning of this section. 

Evidence for  Conjecture~\ref{conj:antidir}  was given in~\cite{ABHLSTR}, where it is verified for all anti\-directed trees of diameter at most~$3$.  In~\cite{ABHLSTR}, it is also noted that in Theorem~\ref{thm:4kn} one can replace the antidirected $k$-edge tree with any antidirected tree whose largest partition class has at most $k$ vertices. 
  In particular,  every $n$-vertex digraph $D$ with more than $2kn$ edges contains each balanced antidirected $k$-edge tree.
  For the antidirected path, 
Klimo\u sov\'a and the author recently improved  the bound on the number of edges of $D$ to roughly $3kn/2$. 

 \begin{theorem}[Klimo\u sov\'a and Stein 2023~\cite{tereza}]\label{antiii}
For each $k\ge 3$, every oriented $n$-vertex graph having at least $(3k-4)n/2$ edges contains each  antidirected path of length $k$. 
 \end{theorem}

 Also,  an approximate version of Conjecture~\ref{conj:antidir} for large dense oriented host graphs and balanced antidirected trees of bounded maximum degree holds.

\begin{theorem}[Stein and Z\'arate-Guer\'en~\cite{camila}]
\label{thm:ES}
For all $\eta\in(0,1)$ and $c\in\mathbb N$, there is a number $n_0\in\mathbb N$  such that for every $n\ge n_0$ and for every $k\ge \eta  n$, every oriented $n$-vertex graph  with more than $(1 + \eta)(k-1)n$ edges contains each balanced antidirected tree $T$ with $k$ edges and $\Delta(T)\leq (\log(n))^c$.
\end{theorem}

 \section{Minimum semidegree}\label{sec:min}\label{sec:semi}

  \subsection{Minimum semidegree and oriented paths in oriented graphs}\label{sec:min-paths}\label{sec:semi1}
  
Many   results in extremal graph theory connect bounds on the minimum degree with the existence of certain subgraphs. One prominent example is  Dirac's theorem from 1952~\cite{dirac}, which states that any graph $G$ on $n\ge 3 $ vertices and of minimum degree at least $n/2$ contains a Hamilton cycle. If we replace the Hamilton cycle with a Hamilton path, i.e.~a path with $n-1$ edges, then the condition on the minimum degree of $G$ can be lowered to $(n-1)/2$. For shorter paths, Dirac, and independently Erd\H os and Gallai (see~\cite{ErdGall59}) observed the following.

\begin{proposition}~\label{factPath}
Let $k\in\mathbb N$. If a graph $G$ has minimum degree at least $k/2$ and a connected component on at least $k+1$ vertices then $G$ contains the path with $k$ edges. 
\end{proposition}
Proposition~\ref{factPath} can be proved with an argument very similar to one of the standard proofs of Dirac's theorem, taking a longest path, and noting that the vertices of this path actually span a cycle, which can then be used to make the path longer.
Observe that the bound on the minimum degree in Proposition~\ref{factPath} is a factor of $1/2$ below the bound needed for a greedy embedding argument.

If we wish to consider the same problem for oriented graphs $D$, we can replace the  minimum degree with the minimum semidegree $\delta^0(D)$.
Note that if $D$ is any oriented graph of minimum semidegree $\delta^0(D) \ge \frac k2$, then the vertices of $D$ have total degree at least $k$, and therefore, each of the components of the underlying undirected graph of~$D$  has more than~$k$ vertices. This phenomenon might be interpreted as an indication that the extra condition in Lemma~\ref{factPath} on the graph having a large component is not necessary for oriented graphs (although for arbitrary digraphs  it may still be needed). 

Jackson \cite{jackson} 
showed a variant of Proposition~\ref{factPath} for directed paths in oriented graphs:

\begin{theorem}[Jackson 1981, Corollary 3 in~\cite{jackson}]\label{thm:jackson}
Let $k\in\mathbb N$. Every oriented graph~$D$ with $\delta^0(D) \ge k$ contains the directed $2k$-edge path.
\end{theorem}
Jackson \cite{jackson} remarks that  this bound on the minimum semidegree is best possible, because of the existence of regular tournaments on $2k+1$ vertices. He also obtains slightly better bounds for strongly connected tournaments.

Note that if we wish to find an odd directed path, say of length $2k-1$, Theorem~\ref{thm:jackson} ensures that a minimum semidegree of at least $k$ is sufficient. In other words, for a directed path of length $k$ (of either parity), a minimum degree exceeding $k/2$ will always be sufficient.

In this spirit, the author of this survey conjectured in~\cite{tree-survey} that the following variant of Jackson's result  holds for any orientation of the path.

\begin{conjecture}[Stein~\cite{tree-survey}] 
\label{conjp}
Let $k\in\mathbb N$. Every oriented graph with  $\delta^0(D) > k/2$ contains every oriented path with $k$ edges.
\end{conjecture}


The bound on the minimum semidegree in Conjecture~\ref{conjp} is best possible. To see this, we can consider, for odd $k$, a $(k-1)/2$-regular tournament on $k$ vertices which does not contain any oriented path with $k$ edges. For even $k$, there are $k$-vertex tournaments of minimum semidegree $k/2-1$. 

A different type of example works for even $k$
and antidirected paths. Consider the {$k/2$-blow-up} of a directed cycle of  length~$\ell\ge 3$, where each vertex $v$ of $C_\ell$ is replaced by an independent set $S_v$ of size $k/2$, and $S_v, S_w$ span a complete bipartite graph with all edges directed from $S_u$ to $S_v$ whenever $vw\in E(C_{\ell})$. Any largest antidirected path in this graph has $k$ vertices, and thus, length $k-1$, while the minimum semidegree of the graph is $k/2$. Note that for this example  it is necessary that $n=\ell k/2$ for some $\ell\ge 3$, in particular, we need $k\le 2n/3$.

Note that Conjecture~\ref{conjp}   becomes very easy if the bound on the minimum semi\-degree is replaced with  $\delta^0(D) \ge k$, as then we can embed any oriented path using a greedy strategy. 
For antidirected paths, Klimo\u sov\'a and the author showed in 2023~\cite{tereza} that 
for each $k\in\mathbb N$ with $k\neq 2$, every 
oriented graph $D$ with $\delta^0(D)\ge \frac{3k-2}{4}$  contains each antidirected path with $k$ edges.
We remark that the case $k=2$ is excluded from Theorem~\ref{thm:KS}, because the bound $\bar \delta^0(D)\ge (6-2)/4=1$ is below the bound from Conjecture~\ref{conjp} and is not sufficient to guarantee an antipath of length two in $D$, since $D$ could be a directed cycle. 

Moreover, an asymptotic version of Conjecture~\ref{conjp} is true for antidirected paths whose length is linear in the order of the host digraph. Namely, 
for all $\eta\in(0,1)$  there is $n_0$ such that for all $n\geq n_0$ and $k\geq \eta n$ every oriented graph $D$ on $n$ vertices with $\delta^0 (D)>(1 + \eta)\frac k2$ contains  every  antidirected path with $k$ edges~\cite{camila}.
This is a direct consequence of 
Theorem~\ref{teo:cam} below.

Some more evidence for  Conjecture~\ref{conjp} can be deduced from results
on semidegree conditions for oriented cycles  in oriented graphs. It was shown by Keevash, Kühn and Osthus \cite{keeko}  that every sufficiently large $n$-vertex oriented graph $G$ of minimum semidegree $\delta^0(G)\geq \frac{3n}8-4$ contains  a  directed Hamilton cycle. Moreover, 
Kelly \cite{kelly} 
proved 
  that every  oriented graph~$D$ of minimum semidegree $\delta^0(D)\geq \frac{3n}8+o(n)$ contains every orientation of a Hamilton cycle (which is tight by examples of H\"aggkvist, see~\cite{kelly}, and has been extended to a pancyclicity result in~\cite{kko}). With every oriented Hamilton cycle, $D$ also contains every oriented spanning path, implying that 
  if $k$ is relatively close to the order of the host, then Conjecture~\ref{conjp} is true, although probably not best possible.
  \begin{question}\label{q:semik}
Let $k,n\in\mathbb N$ with $2n/3<k<n$. What  is the smallest function $f(k)$ such that every $n$-vertex oriented graph with  $\delta^0(D) > f(k)$ contains every oriented path with $k$ edges?
\end{question}
We only ask Question~\ref{q:semik} for $k>2n/3$ because for smaller values of $k$, if $k/2$ divides~$n$, there is the example of the blow-up of the directed $2n/k$-cycle, which fails to contain any antidirected path with $k$ edges.



  \subsection{Minimum semidegree and oriented paths in digraphs}\label{sec:min-paths-di}\label{sec:semi2}

Let us start this section with the known results for semidegrees and Hamilton cycles in digraphs.
In 1960, Ghoulia-Houri~\cite{g-h} proved that a minimum semidegree of at least $n/2$ suffices to guarantee a directed Hamilton cycle in any $n$-vertex digraph~$D$. DeBiasio,  K\"{u}hn,  Molla, Osthus and Taylor~\cite{HamilOrient} showed in 2015 that that same is true for other orientations of the Hamilton cycle, except for antidirected Hamilton cycles. The threshold for antidirected Hamilton cycles is $n/2+1$, as was shown shortly before by DeBiasio and Molla~\cite{HamilAnti}. All of these results are best possible. 

The results from the previous paragraph imply that every $n$-vertex digraph of minimum semidegree at least $n/2$ contains every oriented Hamilton path, since every oriented path can be completed to an oriented cycle which is not antidirected. On the other hand, an $n$-vertex digraph of minimum semidegree $n/2-1$ could be the disjoint union of two complete digraphs of order $n/2$, and therefore not contain any oriented paths with more  than  $n/2-1$ edges.
Similarly, we can consider the disjoint union of complete digraphs of order $k+1$ to see that for all $k$ and $n$ satisfying the obvious divisibility conditions, a minimum semidegree of $k$ in an $n$-vertex digraph is not enough to ensure any oriented path with more than $k$ edges.

So, for possible extensions of Conjecture~\ref{conjp} to digraphs, it will be necessary to require some other condition in addition to the minimum semidegree condition. For instance, one could ask for a lower bound on the order of the largest component of the underlying graph, or equivalently, require that the underlying graph is connected. 
The next result shows that for directed paths, this approach works. Actually, the lower bound on the minimum semidegree can be replaced by a weaker condition: the average of the minimum outdegree and the minimum indegree.
\begin{theorem}[Bermond,   Germa,   Heydemann and Sotteau 1981~\cite{BGHS81}]
\label{thm:bghs}
Let~$D$ be an $n$-vertex digraph whose underlying graph is connected. Then $D$ has a directed path of length $\min\{n-1, \delta^+(D)+\delta^-(D)\}$.
\end{theorem}

Unfortunately, an analogue of Theorem~\ref{thm:bghs} does not hold for antipaths, not even if we require the underlying graph to be strong and of sufficiently large order. To see this, consider two copies of the complete digraph on $k-2$ vertices, $K_1$ and $K_2$. Add two new vertices $v_1$ and $v_2$, and let $v_i$ dominate $V(K_i)$ and be dominated by $V(K_{3-i})$, for each $i=1,2$. The obtained digraph is strong, has order $2k-2$ and minimum semidegree $k-2$, while its longest antidirected path has only $k-1$ edges.
(Note that for $\ell<k/2$,  this example can easily be modified to an $\ell$-strong digraph on $2k -2\ell$ vertices and of minimum semidegree $k-\ell$ whose longest antidirected path has $k-1$ edges.)

So perhaps a different condition on the host digraph has to be added, if we are looking for oriented paths other than directed paths. For antidirected paths, we propose the following.
\begin{question}
Let $D$ be a digraph with $|V(D)|>k$ and $\delta^0(D)>k/2$ such that each pair of vertices is connected by an antidirected path. Does $D$ have an antidirected path of length $k$?
\end{question}
This would even be interesting if we replace $k/2$ with $f(k)$ for some function $f$ with $k/2<f(k)<k$. Note that for $\delta^0(D)\ge k$ the greedy embedding strategy yields an embedding of any oriented path with $k$ edges.

In a slightly different direction, there have been attempts to use bounds on the minimum semidegree to find connectivity-preserving directed paths in digraphs. The next conjecture would be an analogue of   results of Mader for graphs. 
A digraph $D$ is called \defword{$k$-connected} if $|D| \ge k + 1$ and for every pair $a$, $b$ of distinct
vertices, there are at least $k$ internallly disjoint directed paths starting at $a$ and ending in $b$.
\begin{conjecture}[Mader 2012~\cite{maderConTrees}]
\label{conn}
For each $m\in\mathbb N^+$, every $k$-connected digraph~$D$ with $\delta^0(D)\ge 2k+m-1$  contains a $k$-edge directed path $P$ such that deleting $V(P)$ from $D$ leaves a $k$-connected digraph.  
\end{conjecture}
For progress on this conjecture, and an overview of corresponding results in graphs, see~\cite{TLXM}.

  \subsection{Minimum semidegree and oriented trees}\label{sec:min-trees}\label{sec:semi3}
  
In this section, we will try to find oriented trees, other than paths, in oriented graphs and in arbitrary digraphs, employing bounds on the minimum semidegree, and when necessary, additional conditions.

First of all, we observe that the greedy embedding strategy  works for finding any $k$-edge oriented tree, showing that every digraph of  minimum semidegree  at least $k$ contains every oriented tree with $k$ edges. This is no longer true if we change the bound on the minimum semidegree to $k-1$, because of the following example: The $(k-1)$-blow-up of the directed triangle has minimum semidegree $k-1$, but no antidirected star with $k$ edges is present. However,  the antidirected star is  very unbalanced, and the example depends on this feature. It seems natural to suspect that   {balanced} oriented trees with $k$ edges already  appear in 
digraph of  minimum semidegree somewhere below~$k$. 
Moreover, for balanced anti\-directed trees $T$ and oriented host graphs  we suspect this bound should be close to the one we conjecture for antidirected paths.

\begin{conjecture}\label{anticonj}
Every oriented $n$-vertex graph $D$ with $\delta^0(D)>k/2$  contains each balanced antidirected $k$-edge tree of  total maximum degree at most $o(n)$. 
\end{conjecture}

A weaker version would be the conjecture  that every oriented  $n$-vertex graph $D$ with $\delta^0(D)>k/2+o(k)$  contains each balanced antidirected $k$-edge tree of  total maximum degree at most  a power of  $\log n$.  
This is  true  if $D$ is large, and $k$ is large compared to the order of $D$, as we will see in the next result.

\begin{theorem}[Stein and Z\'arate-Guer\'en~\cite{camila}]
\label{teo:cam}
For all $\eta\in(0,1)$ and $c\in \mathbb N$ there is $n_0$ such that for all $n\geq n_0$ and $k\geq \eta n$,  every oriented graph $D$ on $n$ vertices with $\delta^0 (D)>(1 + \eta)k/2$ contains  every balanced antidirected tree $T$ with $k$ edges and with  maximum total degree at most  $(\log(n))^c$.
\end{theorem}
Let us  remark that we have no reason to believe that Theorem~\ref{teo:cam}'s bound on the maximum total degree of $T$ is best possible.
We also note that a more general version of Theorem~\ref{teo:cam}, substituting `oriented graph' with `digraph', is not true, even if we consider a higher bound on the semidegree, as $D$ may decompose into components of order one larger than the semidegree. 

Moreover,  in Theorem~\ref{teo:cam} we cannot omit the condition that $T$ is balanced. The $\ell$-blow-up  of the directed triangle, which has semidegree~$\ell$, does not contain any antidirected tree where one bipartition class has more than $\ell$ vertices. Such antidirected trees exist, with any maximum total degree $\Delta\ge 3$ and $\ell=(\Delta -1)k/\Delta$ (for instance, consider adding $\Delta-2$ leaves to every second vertex of an odd-length path and giving the resulting caterpillar an antidirected orientation).

 For $k=n-1$, 
Theorem~\ref{teo:cam}  follows from a recent result by Kathapurkar and Montgomery~\cite{km}. Extending a  previous result of Mycroft and Naia~\cite{MyNa18} where the tree had constant degree, the authors of~\cite{km} show the following. 

\begin{theorem}[Kathapurkar and Montgomery 2022~\cite{km}]\label{thm:km}
For each $\eta > 0$, there are $c > 0$ and $n_0\in\mathbb N$  such that every $n$-vertex directed graph on $n\ge n_0$ vertices and with minimum semidegree at least $(1/2 + \eta)n$ contains a copy of every $n$-vertex oriented tree of maximum total degree  at most $cn/{\log n}$. 
\end{theorem}

Theorem~\ref{thm:km} generalises a well-known theorem of Koml\'os, S\'ark\"ozy and Szemer\'edi~\cite{KSS2001} for graphs. Indeed, the statement of the result from~\cite{KSS2001} is obtained from the statement of Theorem~\ref{thm:km} by omitting the words `directed', `oriented' and `total', and substituting `semidegree' with `degree'. 

The Koml\'os--S\'ark\"ozy--Szemer\'edi theorem~\cite{KSS2001} is tight in the sense that the maximum degree of the tree could not be of lower order (even if the bound on the minimum degree of the host graph is relaxed). This can be seen by considering the random graph with edge probability $0.9$ and a tree obtained from taking a set of $\log n/c$ stars $K_{1,cn/\log n}$ and joining their centres by new edges in any way to form an $n$-vertex tree, where $c$ is some constant. (See~\cite{KSS2001} for details.)
A direct translation of this construction to digraph shows that also in Theorem~\ref{thm:km}, the bound on the maximum degree of the tree is tight.

Csaba,  Levitt,  Nagy-Gy\"orgy, and Szemer\'edi~\cite{CLNS10} showed in 2010 a variant of the Koml\'os--S\'ark\"ozy--Szemer\'edi theorem~\cite{KSS2001}, with a slightly relaxed bound on the minimum degree of the host graph and a stricter bound on the maximum degree of the tree. Namely, they  showed that for each $\Delta\in\mathbb N$, there are $c > 0$ and $n_0\in\mathbb N$  such that every $n$-vertex  graph on $n\ge n_0$ vertices and with minimum degree at least $\frac n2 + c\log n$ contains  every $n$-vertex  tree of maximum degree  at most $\Delta$. They also showed that this bound on the minimum degree of the host graph is best possible (even for~$\Delta=4$). 

A variant of the result from~\cite{CLNS10} for digraphs seems to be missing yet. We propose the following.
\begin{conjecture}
For every $\Delta$ there is a $c$ such that every large enough $n$-vertex digraph of minimum semidegree at least $n/2 + c\log n$ contains  every $n$-vertex oriented  tree of maximum total degree  at most $\Delta$.
\end{conjecture}

So far, all the results and open problems in this section only applied to  oriented trees of bounded degree.
 However, in  graphs, there are also results 
for containment of trees of unbounded maximum degree. In particular, let us focus on $n$-vertex trees of unbounded maximum degree in $n$-vertex host graphs of large minimum degree. Of course, if we wish to find $n$-vertex stars, then in the host graph, a vertex of degree $n-1$  is needed. It turns out that the presence of such a universal vertex, together with a corresponding condition on the minimum degree, suffices to guarantee all $n$-vertex trees. Indeed, Reed and the author proved in 2023~\cite{RS19a, RS19b} that every $n$-vertex graph of minimum degree exceeding $2n/3$ and maximum degree $n-1$ contains each $n$-vertex tree, if $n$ is large enough. 

A result of this type might also hold for digraphs. Call a vertex $v$ of an $n$-vertex digraph {\it universal} if $d^+(v)=d^-(v)=n-1$.
\begin{question}
Is it true that  every large enough $n$-vertex digraph of minimum semidegree exceeding $2n/3$ that has a universal vertex contains  every oriented  tree on $n$ vertices?
\end{question}

The result from~\cite{RS19a, RS19b} is a special case of an earlier conjecture regarding $k$-edge trees in $n$-vertex host graphs obeying certain minimum/maximum degree conditions. In 2020, Havet, Reed, Wood and the author~\cite{2k3} conjectured that for every $k\in\mathbb N$, every  graph of minimum degree exceeding $2(k+1)/3$ and maximum degree $k$ contains each $k$-edge tree. This bound on the minimum degree is best possible~\cite{2k3} and an asymptotic version for trees of bounded degree has been established in~\cite{BPS1}. 

Thus motivated, we  pose the same problem for digraphs. (One could also see Problem~\ref{lastp}  as a version of 
Theorem~\ref{teo:cam} for oriented trees that not necessarily balanced or antidirected.)
\begin{problem}\label{lastp}
Determine the smallest $f(k)$ such  that for every $k\in\mathbb N$, every  oriented graph (or digraph) of minimum semidegree exceeding $f(k)$ that has a vertex $v$ with $\min\{d^+(v), d^-(v)\}\ge k$ contains each oriented $k$-edge tree.
\end{problem}

Clearly, $f(k)\le k$, and easy modifications of the extremal examples from~\cite{2k3} show that $f(k)\ge 2(k+1)/3$.

\section{Minimum pseudo-semidegree}
\label{sec:pseu}

It is well known that in graphs, a lower bound on the edge densitity can be used to deduce a lower bound on the minimum degree of a subgraph. As we will see now, this is no longer true if we use the semidegree notion in digraphs. However, Lemma~\ref{split} below guarantees that a lower bound on the edge densitity implies a lower bound on the 
\defword{minimum pseudo-semidegree} of a subgraph, where the minimum pseudo-semidegree is a natural variant of the minimum pseudo-semidegree and will be defined below. Later, in Lemma~\ref{antiiiilemma}, we see that for embedding antidirected subgraphs, lower bounds on the minimum semidegree have a very similar effect to lower bounds on the minimum pseudo-semidegree. 

  We  define the \defword{minimum pseudo-semidegree 
$\bar\delta^0(D)$}
of a  digraph $D$ with at least one edge as the maximum~$d$ such that for each vertex $v\in V(D)$ we have $d^+(v)\in \{0\}\cup [d,\infty)$ and  $d^-(v)\in \{0\}\cup [d,\infty)$. In other words, the minimum pseudo-semidegree of a non-empty digraph is the maximum $k$ such that each vertex has either outdegree $0$ or outdegree at least $k$, and has either indegree $0$ or indegree at least $k$. 
The minimum pseudo-semidegree 
of an empty (i.e.~edgeless) digraph is~$0$. 
Note that a digraph of positive minimum pseudo-semidegree can have isolated vertices, i.e.~vertices of in- and outdegree $0$, but not all its vertices can be like this.

\subsection{Edge density and minimum pseudo-semidegree}\label{sec:dens-pseudo}

We start by recalling the relation of edge density and minimum degree of a subgraph in the case of undirected graphs. A folklore lemma states that any graph of high edge density has a subgraph of large minimum degree. More precisely, for each $\ell\in\mathbb N$, every $n$-vertex graph with more than $\ell n$ edges has a subgraph of minimum degree at least $(\ell+1)/2$. This subgraph can be found by successively deleting vertices of low degree.

Because of 
 the example of  the complete bipartite graph on sets $A$ and $B$, with all edges directed from $A$ to $B$, such a lemma cannot exist for digraphs if we replace the subgraph of large minimum degree with a  subgraph of large minimum semidegree (or of large minimum outdegree). However, a similar result is true by Lemma~\ref{split}. 
  This lemma also appears in~\cite{tereza, camila}, and is one of the ingredients for the proof of Theorem~\ref{thm:ES} in~\cite{camila}. It provides a subdigraph of large minimum pseudo-semidegree in a dense digraph.

\begin{lemma}\label{split}
Let $\ell\in\mathbb N$. 
If a digraph $D$ has more than $\ell |V(D)|$ edges, then it contains a digraph $D'$ with $\bar\delta^0(D')\ge (\ell+1)/2$.
\end{lemma}

The  proof of this lemma relies on the observation that on average,
 the vertices of $D$ have in-degree greater than  $\ell$. Also, on average they have out-degree  greater than~$\ell$. We  construct an auxiliary bipartite graph~$B$ by dividing each vertex $v\in V(D)$ into two vertices $v_{in}$ and $v_{out}$, and  letting $v_{in}$ ($v_{out}$) be adjacent to all edges ending (starting) at $v$. We then omit all directions on edges. 
As the average degree of $B$ is greater than   $\ell$, the folklore fact mentioned above shows that~$B$ has a non-empty subgraph $B'$ of minimum degree at least $(\ell+1)/2$. Translating $B'$ back to the digraph setting, it follows that $D$ has a subdigraph~$D'$ with minimum pseudo-semidegree at least $(\ell+1)/2$.  This proves Lemma~\ref{split}.

 \subsection{Minimum semidegree and minimum pseudo-semidegree}
The next  result, Lemma~\ref{antiiiilemma} below, states that the minimum semidegree and the minimum pseudo-semidegree are practically equivalent for the purpose of finding an antidirected subgraph $A$ in an oriented graph $D$, if there is some control over the placement of $A$ in~$D$. 
For digraphs $A$ and $D$, we write $A\subseteq_\gamma D$
if for each set $V^*\subseteq V(D)$ of size at least $\gamma |V(D)|$ and for each $x\in V(A)$, there is an embedding of $A$ in $D$ with $x$ mapped to $V^*$.

Part~(i) of the following lemma appears as Lemma 7.2 in~\cite{camila}, and part (ii) can be proved in the same way.

\begin{lemma}\label{antiiiilemma} For any antidirected graph $A$ whose underlying graph is connected and for any
$\ell, n_0\in \mathbb N$ the following holds.
\begin{enumerate}[(i)]
\item
If for each oriented graph $D$ on at least $n_0$ vertices with $\delta^0(D)\ge \ell$ we have $A\subseteq_{1/8} D$, then  each oriented graph $D'$ on at least~$n_0$ vertices with $\bar\delta^0(D')\ge \ell$ contains $A$. 
\item 
If for each digraph $D$ on at least $n_0$ vertices with $\delta^0(D)\ge \ell$ we have $A\subseteq_{1/8} D$, then  each digraph $D'$ on at least~$n_0$ vertices with $\bar\delta^0(D')\ge \ell$ contains $A$. 
\end{enumerate}
\end{lemma}

The idea of the proof of Lemma~\ref{antiiiilemma} is as follows. Given $D'$ and $A$, we can assume that $D'$ has no isolated vertices, and we will construct a graph $D$ by gluing together four copies of $D'$, two of them with reversed edge directions. 

More precisely, we take a copy of $D'$ as it is and another copy of $D'$ with all edge directions reversed. We identify the vertices of first copy that have outdegree $0$ with the corresponding vertices in the second digraph (in that graph, these vertices have indegree $0$). This gives a digraph $D''$ that still has  minimum pseudo-semidegree $\ell$. We now repeat the procedure from above, taking a second copy of $D''$ with all edge direction reversed, and identifying the vertices of outdegree $0$ from the first copy with their copies in the second digraph.  We also identify  the vertices of indegree $0$ from the first copy with their copies in the second digraph. The obtained digraph $D$ has minimum semidegree $\ell$. As $A\subseteq_{1/8} D$, it is not hard to see that we can find $A$ in $D$, guaranteeing that it is fully contained in one of the `properly oriented' copies of $D'$. For details, see~\cite{camila}.

\subsection{\hskip-.35cm Minimum pseudo-semidegree and oriented paths and trees}

The results of the previous subsection, in particular Lemma~\ref{antiiiilemma}, hint at the possibility that if we are looking for antidirected paths or trees, then instead of a lower bound on the minimum semidegree of the host digraph, as in Section~\ref{sec:semi}, we could work with a lower bound on the minimum pseudo-semidegree of the host digraph, which would be a weaker condition.

So, recalling Conjecture~\ref{conjp} from Section~\ref{sec:semi}, the following conjecture seems plausible.
\begin{conjecture}\label{pseudoanti}
For each $k\in\mathbb N$, every oriented graph with  $\bar\delta^0(D) > k/2$ contains every antidirected path with $k$ edges.
\end{conjecture}
Observe that Conjecture~\ref{pseudoanti} is not true for other orientations of the path, because of the example of the complete bipartite graph $(A,B)$ with all edges directed from~$A$ to $B$. Also note that replacing the bound  of $k/2$ in Conjecture~\ref{pseudoanti} with $k$ makes the conjecture trivial, because of the greedy embedding argument. 

As evidence for Conjecture~\ref{pseudoanti}, we now present analogues of some of the results from Section~\ref{sec:semi}. The first one of these uses a bound on the minimum pseudo-semidegree that is  close to $3k/4$, which lies halfway betweeen the trivial bound and the bound from the conjecture.

\begin{theorem}[Klimo\u sov\'a and Stein 2023~\cite{tereza}]\label{thm:KS}
Let $k\in\mathbb N$ with $k\neq 2$. Every 
oriented graph $D$ with $\bar\delta^0(D)> \frac{3k-2}{4}$  contains each antidirected path with $k$ edges.
\end{theorem}

We remark that  one can use Lemma~\ref{split} to see that Theorem~\ref{thm:KS} implies Theorem~\ref{antiii}.

Moreover,  for all $\eta\in(0,1)$  there is $n_0$ such that for all $n\geq n_0$ and $k\geq \eta n$ every oriented graph $D$ on $n$ vertices with $\bar\delta^0 (D)>(1 + \eta)\frac k2$ contains  every  antidirected path with $k$ edges. This  follows directly from Theorem~\ref{pseudocam}, which we state next.

\begin{theorem}
\label{pseudocam}
For all $\eta\in(0,1)$ and $c\in \mathbb N$ there is $n_0$ such that for all $n\geq n_0$ and $k\geq \eta n$,  every oriented graph $D$ on $n$ vertices with $\bar\delta^0 (D)>(1 + \eta)k/2$ contains  every balanced antidirected tree $T$ with $k$ edges and with  maximum total degree at most  $(\log(n))^c$.
\end{theorem}

Theorem~\ref{pseudocam} follows immediately from Theorem 5.2 in~\cite{camila} and from  Lemma~\ref{antiiiilemma} above.

We believe that the approximation and the condition that $k$ and $n$ are large can be omitted from Theorem~\ref{pseudocam}:
\begin{conjecture}\label{pseudoantitree}
For each $k\in\mathbb N$, every oriented graph on $n$ vertices with  $\bar\delta^0(D) > k/2$ contains every antidirected balanced tree having $k$ edges and  maximum total degree~$o(n)$.
\end{conjecture}
%
%
%
%
%
%
%
%


%
%
%
%



 \bibliographystyle{amsplain}
 \bibliography{trees2.bib}



\end{document}